\input amstex
\documentstyle{amsppt}
\topmatter
\title Hahn-Banach operators
\endtitle
\rightheadtext{Hahn-Banach operators}

\author M.~I.~Ostrovskii
\endauthor

\address  Department of Mathematics,
The Catholic University of America, 
Washington, D.C. 20064, USA
\endaddress

\email ostrovskii\@cua.edu\endemail

\abstract We consider real spaces only. 

{\bf Definition.} An operator $T:X\to Y$ between Banach spaces
$X$ and $Y$ is called a
{\it Hahn-Banach operator}
if for every isometric embedding of the space $X$
into a Banach space $Z$ there exists a norm-preserving
extension $\tilde T$ of $T$ to $Z$.

A geometric property of Hahn-Banach operators of finite rank acting
between finite-dimensional normed spaces is found. This property is
used to characterize pairs of finite-dimensional normed spaces
$(X,Y)$ such that there exists a Hahn-Banach operator $T:X\to Y$
of rank $k$. The latter result is a generalization of a recent result
due to B.~L.~Chalmers and B.~Shekhtman.
\endabstract
\endtopmatter
\document

Everywhere in this paper we consider only real linear spaces.
Our starting point is the classical Hahn-Banach theorem
(\cite{H}, \cite{B1}). The form of the Hahn-Banach theorem
we are interested in can be stated in the following way.

\proclaim{Hahn-Banach Theorem} Let $X$ and $Y$ be Banach spaces,
$T:X\to Y$ be a bounded linear operator of rank $1$ and $Z$ be a
Banach space containing $X$ as a subspace. Then there exists a
bounded linear operator $\tilde T:Z\to Y$ satisfying
\smallskip

{\rm (a)} $||\tilde T||=||T||$;
\smallskip

{\rm (b)} $\tilde Tx=Tx$ for every $x\in X$.
\endproclaim

\definition{Definition 1} An operator $\tilde T: Z\to Y$ satisfying (a)
and (b) for a bounded linear operator $T:X\to Y$
is called a {\it norm-preserving extension} of $T$ to $Z$.
\enddefinition

The Hahn-Banach theorem is one of the basic principles
of linear analysis. It is quite
natural that there exists a vast literature on generalizations of the
Hahn-Banach theorem for operators
of higher rank. See the papers G.~P.~Akilov \cite{A},
J.~M.~Borwein \cite{Bor}, B.~L.~Chalmers
and B.~Shekhtman \cite{CS}, G.~Elliott and I.~Halperin \cite{EH},
D.~B.~Goodner \cite{Go}, A.~D.~Ioffe \cite{I}, S.~Kakutani \cite{Kak},
J.~L.~Kelley \cite{Kel}, J.~Linden\-strauss \cite{L1}, \cite{L2},
L.~Nachbin \cite{N1}
and M.~I.~Ostrov\-skii \cite{O}, representing different directions of
such generalizations, and references therein.
There exist two interesting surveys devoted to the
Hahn-Banach theorem and its generalizations, see G.~Buskes \cite{Bus} and
L.~Nachbin \cite{N2}.
\medskip

We shall use the following natural definition.

\definition{Definition 2} An operator $T:X\to Y$ between Banach spaces
$X$ and $Y$ is called a
{\it Hahn-Banach operator}
if for every isometric embedding of the space $X$
into a Banach space $Z$ there exists a norm-preserving
extension $\tilde T$ of $T$ to $Z$.
\enddefinition

The existence of non-Hahn-Banach operators was mentioned in the remarks
to Chapter IV of Banach's book, see \cite{B2, p.~234}.
S.~Banach and S.~Mazur \cite{BM} proved that the identity operator on $l_1$
is a non-Hahn-Banach operator (in fact,
this operator does not have even continuous extensions for some
isometric embeddings). It has been
known for a long time 
that there exist non-Hahn-Banach operators of rank $2$
(see F.~Bohnenblust \cite{Boh}, an important relevant result
was proved relatively recently by  H.~K\"onig \&\
N.~Tomczak-Jaegermann \cite{KT}). A problem of characterization of
Hahn-Banach operators arises in a natural way.
\medskip

Factorizational characterizations of Hahn-Banach operators
are well known. In particular,
using by now a standard technique (that goes back to G.~P.~Akilov
\cite{A}, D.~B.~Go\-odner \cite{Go, pp.~92--93} and R.~Phillips \cite{P,
p.~538})
it is easy to show
that an operator $T:X\to Y$ is a Hahn-Banach operator if and only if
for some set $\Gamma$ 
there exist operators $T_1:X\to l_\infty(\Gamma)$
and $T_2:l_\infty(\Gamma)\to Y$ such that
$T_2T_1=T$ and $||T_2||||T_1||=||T||$. (See \cite{J} for
the undefined terminology from the theory of Banach spaces.)
\medskip

One of the main purposes of the present paper is to find a
geometric property of Hahn-Banach operators of finite rank
acting between finite-dimensional normed spaces (see Theorem 1).
This property does not imply that the operator is a Hahn-Banach
operator (see the remark after Theorem 1), but it can be used to
answer the following
question: given $k\in\Bbb N$, for which pairs of finite-dimensional spaces
$(X,Y)$ does there exist a Hahn-Banach operator $T:X\to Y$ of rank $k$?
(See Theorem 2.) This result is a generalization of a
recent result due to B.~L.~Chalmers and B.~Shekhtman
\cite{CS}.
\medskip

\remark{Remark} Let $\alpha\in\Bbb R,\ \alpha\ne 0$. It is clear that
$\alpha T$ is a Hahn-Banach operator if and only if $T$ is a Hahn-Banach
operator. Hence studying Hahn-Banach operators it is enough to consider
Hahn-Banach operators of norm $1$.
\endremark
\medskip

We need the following notation. By $S(X)$ and $B(X)$ we denote
the unit sphere and the unit ball of a Banach space $X$ respectively.
Let $X$ be a finite-dimensional Banach space. An intersection of $B(X)$ with
a supporting  hyperplane of $B(X)$ will be called
a {\it support set} of $B(X)$. By the {\it dimension} of a set in a
finite-dimensional space we mean the dimension of its affine hull.
(See \cite{S} for the undefined terminology from the theory
of convex bodies.)
We define $f(X)$ to be the maximal dimension of 
a support set of $B(X)$. For $x\in S(X)$ we define $d(x)$ to
be the dimension of the set  $\{x^*\in S(X^*):\ x^*(x)=1\}$.
It is clear that $d(x)=0$ if and only if $x$ is a smooth point;
in the general case $d(x)$ indicates the number of linearly
independent directions of non-smoothness of the norm at $x$.

\proclaim{Theorem 1} Let $X$ and $Y$ be finite-dimensional Banach spaces
and $T:X\to Y$ be a Hahn-Banach operator of rank $k$. Assume that
$||T||=1$ and let $x_0\in S(X)$ be such that $||Tx_0||=1$. Then
$Tx_0$ belongs to a support set of $B(Y)$ of dimension  $\ge k-1-d(x_0)$.
\endproclaim

\demo{Proof}
Let $C(S(X^*))$ denote the space of all continuous functions
on $S(X^*)$ with the $\sup$ norm. 
We identify $X$ with a subspace of
$C(S(X^*))$ in the following way: every vector
is identified with its restriction  (as a function
on $X^*$) to $S(X^*)$. We introduce the following notation:
$C=C(S(X^*))$ and $B_C=B(C(S(X^*)))$.
\medskip

Since $T$ is a Hahn-Banach operator, there exists
$\tilde T:C\to Y$ such that $\tilde T|_X=T$ and $||\tilde T||=1$.
We shall use $\tilde T$ to find a ``large'' support set of
$B(Y)$.
\bigskip

Since $||Tx_0||=1$, there exists $h\in S(Y^*)$ such
that $h(Tx_0)=1$. Let $F=\{x^*\in S(X^*):\ x^*(x_0)=1\}$.
Observe that
$$T^*h\in F.\eqno{(1)}$$ 
\medskip

Choose a basis $\{y_1,\dots,y_m\}$ in $Y$ such that
$y_1=Tx_0$ and
$y_2,\dots,y_m\in\ker h.$
The operator $\tilde T$ can be represented in the form
$$\tilde T=\sum_{i=1}^m\mu_i\otimes y_i,\ \ \mu_i\in C^*.$$
By the F.Riesz representation theorem (see e.g. \cite{DS},
p.~265) we may identify
$\mu_i$ with (signed) measures on $S(X^*)$.
\medskip

Our first purpose is to show that $\mu_1$ is supported
on $F\cup(-F)$.
\medskip

We have $\tilde T(B_C)\subset B(Y)\subset
\{y:\ |h(y)|\le 1\},\ h(y_1)=1$ and $y_2,y_3,\dots, y_m\in \ker h$.
Therefore for every $z\in B_C$ we have
$$\mu_1(z)=h(\sum_{i=1}^m\mu_i(z)y_i)=h(\tilde T(z))$$
and $|\mu_1(z)|=|h(\tilde T(z))|\le 1$. Hence
$$||\mu_1||\le 1\eqno{(2)}$$
Also, since $\tilde Tx_0=y_1$, we have
$$\mu_1(x_0)=1.\eqno{(3)}$$

Conditions (2), (3) and $||x_0||=1$ imply that
$\mu_1$ is supported on $F\cup(-F)$. (By this we mean
that the restriction of $\mu_1$ to
$S(X^*)\backslash(F\cup(-F))$ is a zero measure.)
\medskip

We decompose 
$\mu_i=\nu_i+\omega_i$, where $\nu_i$ is the restriction
of $\mu_i$ to $F\cup(-F)$. Since
$\mu_1$ is supported on $F\cup(-F)$, then $\omega_1=0$.
\medskip

Since $T$ is of rank $k$, there exists a subspace
$L\subset X$ of dimension $k$ such that $T|_L$
is an isomorphism. Let
$$M=\{x\in L:\ 
\forall x^*\in F,\ \ x^*(x)=0\}.$$
Then $\dim M\ge k-d(x_0)-1$.
\medskip

Let $x\in B(M)$. The definitions of $M$ and $\nu_i$ imply that
$$\nu_i(x)=0,\ i\in\{1,\dots,m\}.\eqno{(4)}$$
(Recall that we identify vectors in $M$ with the corresponding
functions in $C$.)
\smallskip

Now we construct a ``mixture'' of $x$ and $x_0$.
\smallskip

It is clear that for each $\delta>0$ there exists a function
$g_\delta\in B_C$ such that 
$$g_\delta|_{F\cup(-F)}=x_0|_{F\cup(-F)}$$
and the
restrictions of $g_\delta$ and $x$ to the complement of the $\delta-$neighbourhood of $F\cup(-F)$ coincide.
\medskip

We have
$$\lim_{\delta\downarrow 0}\tilde Tg_\delta=\lim_{\delta\downarrow 0}
\sum_{i=1}^m\mu_i(g_\delta)y_i=
\lim_{\delta\downarrow 0}\sum_{i=1}^m(\nu_i(g_\delta)+
\omega_i(g_\delta))y_i.$$

We have 
$\nu_i(g_\delta)=\nu_i(x_0)$ for every $\delta>0$ and $i\in\{1,\dots, m\}$.

It is clear that $\omega_i(F\cup(-F))=0$. By the definition
of $g_\delta$ it follows that
$\lim_{\delta\downarrow 0} g_\delta(x^*)=x(x^*)$ for $x^*\in S(X^*)
\backslash(F\cup(-F))$ and that the functions
$g_\delta$ are uniformly bounded. By the
Lebesgue dominated convergence theorem we get
$$\lim_{\delta\downarrow 0}\omega_i(g_\delta)=\omega_i(x).$$

Therefore

$$\lim_{\delta\downarrow 0}\tilde Tg_\delta=
\sum_{i=1}^m(\nu_i(x_0)+\omega_i(x))y_i.$$

Equation (4) implies that $\nu_i(x_0)=\nu_i(x_0+x)$.
Using this and the fact that $\omega_1=0$, we get
$$\lim_{\delta\downarrow 0}\tilde Tg_\delta=
\sum_{i=1}^m\mu_i(x_0+x)y_i-\sum_{i=2}^m\omega_i(x_0)y_i=
T(x_0+x)-\sum_{i=2}^m\omega_i(x_0)y_i.$$
\bigskip

Since $g_\delta\in B_C$, $||\tilde T||=1$ and $B(Y)$ is closed, then
$$T(x_0+x)-\sum_{i=2}^m\omega_i(x_0)y_i\in B(Y)$$
for every $x\in B(M)$.

By (1) we have $h(Tx)=0$ for every $x\in M$. Recall, also, that
$y_2,\dots,y_m\in\ker h$. Therefore
$$h\left(T(x_0+x)-\sum_{i=2}^m\omega_i(x_0)y_i\right)=
hTx_0=1$$
for every $x\in M$.
Since $T|_M$ is an isomorphism and the vector
$\sum_{i=2}^m\omega_i(x_0)y_i$
does not depend on $x$, the intersection of 
$B(Y)$ with the supporting hyperplane
$\{y:\ h(y)=1\}$ has dimension $\ge\dim M\ge k-d(x_0)-1$.
$\square$
\enddemo

\proclaim{Corollary} The existence of a Hahn-Banach operator
$T:X\to Y$ of rank $k$ implies
$f(X^*)+f(Y)\ge k-1$.
\endproclaim

\demo{Proof} Observe that $f(X^*)\ge d(x_0)$ and that $f(Y)$ is
greater or equal to the dimension of any support set of $B(Y)$
containing $Tx_0$. $\square$
\enddemo

\remark{Remark} There exist operators satisfying the condition
of Theorem 1 that are not Hahn-Banach. In fact, let $T$ be the
identity mapping of $X=l_1^n$ onto the space $Y$ whose unit
ball is the intersection of $(1+\varepsilon)B(l_1^n),\
(\varepsilon>0)$ and $B(l_\infty^n)$. It is easy to see that

(1) the norm of this operator is $1$;

(2) the only points where the operator attains its norm are
$\pm e_1,\pm e_2,\dots,\pm e_n$, where $\{e_1,\dots,e_n\}$ is
the unit vector basis;

(3) the points $\pm e_1,\pm e_2,\dots,\pm e_n$ are contained
in $(n-1)$-dimensional support sets of $B(Y)$.
\smallskip

Therefore $T$ satisfies the condition of Theorem 1 with $k=n$.
On the other hand, the operator $T$ is not a Hahn-Banach operator
if $n\ge 3$ and $\varepsilon$ is small enough. In fact, 
$||T^{-1}||=1+\varepsilon$. Therefore, if $T$ were a Hahn-Banach operator
it would imply that for every Banach space $Z$
containing $l_1^n$ as a subspace there exists a projection
onto $l_1^n$ with the norm  $\le 1+\varepsilon$. 
It remains to apply the well-known result of B.~Gr\"unbaum
(see \cite{Gr} or \cite{J, p.~81}).
\endremark
\medskip

B.~L.~Chalmers and B.~Shekhtman \cite{CS} characterized 2-dimensional
spaces having isomorphisms that are Hahn-Banach operators. One of the
steps in their approach is an embedding of the considered spaces into
$L_1$. This is why they got the restriction on the dimension
(spaces of dimension $\ge 3$ may be non-isometric to any subspace of
$L_1$, see J.~Lindenstrauss \cite{L2, p.~494}).
Our next purpose is to extend their results and to characterize
pairs ($X,Y$) of finite-dimensional normed linear
spaces such that there exists a Hahn-Banach operator
$T:X\to Y$ of rank $k$.
\bigskip

\proclaim{Theorem 2} Let $X$ and $Y$ be finite-dimensional
normed linear spaces and let $k$ be a positive integer satisfying
$k\le\min\{\dim X, \dim Y\}$. There exists a Hahn-Banach operator
$T:X\to Y$ of rank $k$ if and only if
$f(X^*)+f(Y)\ge k-1$.
\endproclaim

\demo{Proof} The necessity has been already proved
(see the corollary).
\bigskip

Sufficiency. Suppose that
$$k\le\min\{\dim X,\dim Y, f(X^*)+f(Y)+1\}.$$

It is well known (see e.g. \cite{KS}) that in order to show that $T:X\to Y$
is a Hahn-Banach operator it is enough to show that it has a norm-preserving
extension to the space $C=C(S(X^*))$ (The space $X$
is embedded into $C$ in the same way as in Theorem 1).
Therefore, if an operator $Q:C\to Y$ is such that the
restriction of $Q$ to $X$ has rank $k$ and
$||Q||=||Q|_X||=1,$
then $T=Q|_X$ is a Hahn-Banach operator of rank $k$.
\medskip

Our purpose is to construct such $Q$.
Let $n=f(Y)$ and $m=f(X^*)$. Let $y_0,y_1,\dots,y_n\in Y$ be linearly
independent and such that
$$\{y:\ y=\theta y_0+\sum_{i=1}^na_iy_i,\hbox{ where }
\theta=\pm 1,\ |a_i|\le 1\}\subset S(Y).$$

Let $x^*_0,x^*_1,\dots,x^*_m\in X^*$ be linearly
independent and such that
$$\{x^*:\ x^*=\theta x^*_0+\sum_{i=1}^mb_ix^*_i,\hbox{ where }
\theta=\pm 1,\ |b_i|\le 1\}\subset S(X^*).$$

Let $x_0\in S(X)$ be such that $x^*_0(x_0)=1$. Let
$$x^*_0, x^*_1,\dots,x^*_m,x^*_{m+1},\dots, x^*_r,$$
where $r=\dim X-1$ be a basis in $X^*$ satisfying the condition
$x^*_{m+1}(x_0)=\dots=x^*_r(x_0)=0$. (Observe that the condition
$x^*_1(x_0)=\dots=x^*_m(x_0)=0$ follows from our choice of the vectors.)
Let $x_0,x_1,\dots,x_r$ be its biorthogonal vectors.
\bigskip

Let $y_0,y_1,\dots,y_s$, where $s=\dim Y-1$, be a
basis in $Y$.

We suppose that $k>m+1$. (It will be clear from our argument
which changes should be made if it is not the case.)
\medskip

We define an operator $Q_1:C\to Y$ as follows. Let
$\mu_0, \mu_{m+1},\dots,\mu_k$ be norm-preserving extensions
of $x^*_0, x^*_{m+1}, \dots, x^*_k$ to $C$. Let
$$Q_1(f)=\mu_0(f)y_0+
\sum_{i=1}^{k-m-1}\frac{\mu_{m+i}(f)}{||\mu_{m+i}||}y_i.$$

It is clear that for $x\in X$
$$Q_1(x)=x^*_0(x)y_0+
\sum_{i=1}^{k-m-1}\frac1{||\mu_{m+i}||}x^*_{m+i}(x)y_i.\eqno{(5)}$$

We have supposed that $k-m-1\le n$. This and the
choice of $y_0,\dots,y_n$ implies that $||Q_1||\le 1$.
\bigskip

Our next step is to show that
there exist signed measures $\nu_0,\nu_1,\dots, \nu_m$
of norm 1 on $S(X^*)$ satisfying the conditions
$$\nu_i(x_j)=\delta_{i,j},\ i=0,\dots,m,\ j=0,\dots,r
\eqno{(6)}$$
and
$$\forall f\in B_C\  \forall j\in\{1,\dots, m\}\ |\nu_j(f)|\le
1-|\nu_0(f)|.\eqno{(7)}$$

Let us verify that the following measures satisfy these
conditions.
We define $\nu_j,\ j=1,\dots, m$ as atomic measures with atoms at
$x_0^*+\sum_{i=1}^m\theta_ix^*_i,\ \theta_i=\pm 1$
satisfying $\nu_j(x_0^*+\sum_{i=1}^m\theta_ix^*_i)=2^{-m}\theta_j$
and $\nu_0$ as an atomic measure satisfying 
$\nu_0(x_0^*+\sum_{i=1}^m\theta_ix^*_i)=2^{-m}$.
\medskip

Condition (6) follows from the fact that the sequences
$\{x_0,\dots,x_r\}$ and\newline $\{x_0^*,\dots,x_r^*\}$ are
biorthogonal.
\medskip

Let us verify condition (7).
Denote by $\Bbb I$ the function that is identically $1$ on
$S(X^*)$. Let $f\in B_C$, $j\in\{1,\dots,m\}$. Then $\nu_j(-f)=
\nu_j(\Bbb I-f)\le$ (since the function $\Bbb I-f$ is nonnegative)
$\le\nu_0(\Bbb I-f)=1-\nu_0(f).$ (Here we explicitly use the fact that
the spaces are real.)
\medskip

This proves (7) in the case when $\nu_j(f)$ is negative
and $\nu_0(f)$ is positive. It remains to observe that
\smallskip

A. Since we may consider $-f$ instead of $f$
it is enough to prove (7) for
functions with positive $\nu_0(f)$.
\smallskip

B. For each function $f\in B_C$ and $j\in\{1,\dots, m\}$
there exists a function $f_j\in B_C$ such that
$\nu_0(f_j)=\nu_0(f)$ and $\nu_j(f_j)=-\nu_j(f)$.
\bigskip

We introduce $Q_2:C\to Y$
by the equality
$$Q_2(f)=\nu_0(f)y_0+\sum_{i=1}^m\alpha_i\nu_i(f)y_{k-m-1+i}.$$

Condition (7) implies that there
exists $\{\alpha_i\}_{i=1}^m$ such that $\alpha_i\ne 0$ for every $i$ and
$||Q_2||\le 1$.
\medskip

Condition (6) implies that for $x\in X$ we have
$$Q_2(x)=x^*_0(x)y_0+\sum_{i=1}^m\alpha_i x^*_i(x)y_{k-m-1+i}.\eqno{(8)}$$

Now, let $Q=\frac12(Q_1+Q_2)$. Our estimates for $||Q_1||$
and $||Q_2||$ immediately imply $||Q||\le 1$. On the other
hand, equations (5) and (8) imply that $Q(x_0)=y_0$.
Hence $||Q|_X||\ge 1$ and $||Q||=||Q|_X||=1$.

Also, from (5) and (8) we get for every $x\in X$:
$$Q(x)=x^*_0(x)y_0+\sum_{i=1}^{k-m-1}\frac1{2||\mu_{m+i}||}x^*_{m+i}(x)y_i
+\sum_{i=1}^m\frac12\alpha_ix^*_i(x)y_{k-m-1+i}.$$

Since the sequences $\{y_0\dots,y_s\}$ and
$\{x_0^*,\dots,x_r^*\}$ are linearly independent, it follows that
$Q|_X$ is of rank $k$. $\square$
\enddemo
\bigskip

{\bf Acknowledgement.} The author would like to thank Prof. Bruce
L. Chalmers for giving a copy of \cite{CS} and for suggesting the
problem considered here.

\enddocument